\documentclass{amsart}

\usepackage[T1]{fontenc}
\usepackage{latexsym,amssymb,amsmath}

 \newcommand{\Galg}{\mathbf{G}}
 
 \newcommand{\Balg}{\mathbf{B}}

 \newcommand{\Ualg}{\mathbf{U}}

 \newcommand{\Talg}{\mathbf{T}}

 \newcommand{\Ind}{\operatorname{Ind}}

 \newcommand{\Nn}{\operatorname{N}}

 \newcommand{\Cen}{\operatorname{C}}
 
 \newcommand{\Sh}{\operatorname{Sh}}

 \newcommand{\semi}[1]{\rtimes\langle\,#1\,\rangle}
 
 \newcommand{\cyc}[1]{\langle\,#1\,\rangle}

 \newcommand{\F}{\mathbb{F}}
 
 \newcommand{\N}{\mathbb{N}}
 
 \newcommand{\Q}{\mathbb{Q}}

 \newcommand{\Irr}{\operatorname{Irr}}

\newtheorem{theorem}{Theorem}[section] 
\newtheorem{lemma}[theorem]{Lemma}     

\newtheorem{proposition}[theorem]{Proposition}

\newtheorem{remark}[theorem]{Remark}

\title[On the unipotent characters of the Ree groups]
 {On the unipotent characters of the Ree groups of type~$G_2$}

\author{Olivier Brunat}

\address{Ruhr-Universit\"at Bochum\\
Fakult\"at f\"ur Mathematik\\
Raum NA 2/33\\
D-44780 Bochum\\}
\email{Olivier.Brunat@ruhr-uni-bochum.de}

\begin{document}
\maketitle
\begin{abstract}
This note is concerned with the unipotent characters of the Ree groups of
type $G_2$. We determine the roots of unity associated 
by Lusztig and Digne-Michel to
unipotent characters of $^2G_2(3^{2n+1})$ and we prove that the
Fourier matrix of $^2G_2(3^{2n+1})$  defined by Geck and Malle 
satisfies
a conjecture of Digne-Michel. Our main tool is the Shintani descent of
Ree groups of type $G_2$.
\end{abstract}

\section{Introduction}\label{intro}

Let $\Galg$ be a connected reductive group defined over the finite field
$\F_q$ with $q$ elements of characteristic $p>0$, and let $F$ be
the corresponding Frobenius map. Let 
$\Talg$ be a maximal rational torus of $\Galg$, contained in a rational Borel
subgroup $\Balg$ of $\Galg$. We denote by
$W=\Nn_{\Galg}(\Talg)/\Talg$ the Weyl group of $\Galg$. 
For $w\in W$, there is the
corresponding Deligne-Lusztig 
character $R_w$ of the finite fixed-point subgroup $\Galg^F$. 
We refer to~\cite[\S7.7]{Carter2} for a
precise construction. Then we define the set of
unipotent characters of $\Galg^F$ by
$$\mathcal{U}(\Galg^F)=\{\chi\in\Irr(\Galg^F)\ |\
\cyc{\chi,R_w}_{\Galg^F}\neq 0\ \textrm{for some }w\in W\}.$$

Lusztig~\cite{LusBook} and Digne-Michel~\cite{DMShin} associated to
every $\chi\in\mathcal{U}(\Galg^F)$ a root of
unity~$\omega_{\chi}\in\overline{\Q}_{\ell}$ (where 
$\overline{\Q}_{\ell}$ is the $\ell$-adic field for $\ell\neq p$) as
follows.
We denote by $\delta$
the order of the automorphism of $W$ induced by $F$.
For $w\in W$, we denote by $n_w\in\Nn_{\Galg}(\Talg)$ an element such
that $n_w\Talg=w$. We set $X_w=\{x\Balg\ |\ x^{-1}F(x)\in \Balg
n_w\Balg\}$ the corresponding Deligne-Lusztig variety.
For every integer $j$ we
denote by $H_c^j(X_w,\overline{\Q}_{\ell})$ the $j$-th $\ell$-adic
cohomology space with compact support over 
$\overline{\Q}_{\ell}$, associated to $X_w$.  
Therefore, the groups $\cyc{F^\delta}$ and $\Galg^F$ act on $X_w$, which induces
linear operations on
$H_c^j(X_w,\overline{\Q}_{\ell})$. Hence
$H_c^j(X_w,\overline{\Q}_{\ell})$ is a
$\overline{\Q}_{\ell}\Galg^F$-module, and $F^\delta$ acts on this space as
a linear endomorphism. We also fix an eigenvalue $\lambda$ of $F^\delta$,
and we denote by $F_{\lambda,j}$ its generalized eigenspace. Since the
actions of $\cyc{F^\delta}$ and $\Galg^F$ commute, the space
$F_{\lambda,j}$ is a $\overline{\Q}_{\ell}\Galg^F$-module. Moreover,
the irreducible constituents of $\Galg^F$ occurring in this module are
unipotent.
Conversely, for every unipotent character $\chi$  of $\Galg^F$, there are $w\in
W$, $\lambda\in\overline{\Q}_{\ell}^{\times}$ and $j\in\N$, such that
$\chi$ occurs in the character associated to $F_{\lambda,j}$.
Lusztig~\cite{LusBook} and Digne-Michel~\cite{DMShin} 
have shown that $\lambda=q^{s/2}\omega_{\chi}$ for some non-negative
integer $s$ and a root of unity
$\omega_{\chi}\in\overline{\Q}_{\ell}^\times$
which depends only on $\chi$. 

The set of roots associated as above to the unipotent characters
have been computed by
Lusztig in~\cite{LusBook} for finite reductive groups if
 $F$ is split. Moreover Lusztig computed
in~\cite{LusCox} the set of roots for the unipotent characters appearing in $
H_c^j(X_{w_{\operatorname{cox}}},\overline{\Q}_{\ell})$ where
$w_{\operatorname{cox}}$ denotes the Coxeter element of $W$, 
with no condition on~$F$. This work is completed for the cases
that $F$ is a
non-split Frobenius map by Geck and Malle in~\cite{GM}.
However, for few cases, the methods of Lusztig and Geck-Malle 
allow to associate the roots to their unipotent characters only up 
to complex conjugation. This is for example the case for the Suzuki and
the Ree groups. In~\cite{Br3} we remove this indetermination for the
unipotent characters of the 
Suzuki groups.

Moreover, we recall that Lusztig~\cite{LusBook} associated to most of sets
$\mathcal{U}(\Galg^F)$ some non-abelian Fourier matrices, which involve the
decomposition in unipotent characters of $R_w$ for $w\in W$.

This note is concerned with the Fourier matrices and 
the roots of unity, associated as above to
the unipotent characters of the Ree group~$^2G_2(q)$ for $q=3^{2n+1}$.
For these groups, the method in~\cite{LusBook} does not allow to define 
Fourier matrices. However using the theory of character sheaves, Geck
and Malle give a more general definition for these
matrices~\cite[5.1]{GM}. For the Ree groups of type $G_2$, they obtained 
the following matrix~\cite[5.4]{GM}
$$
\frac{\sqrt{3}}{6}
\begin{bmatrix}
1&1&1&1&2&2\\
\sqrt{3}&-\sqrt{3}&\sqrt{3}&-\sqrt{3}&0&0\\
0&0&2&2&0&-2\\
2&2&0&0&-2&0\\
1&1&-1&-1&2&-2\\
\sqrt{3}&-\sqrt{3}&-\sqrt{3}&\sqrt{3}&0&0\\
\end{bmatrix}.
$$

We set $I=\{1,3,5,6,7,8,9,10\}$. The Ree group  $^2G_2(q)$ 
has $8$ unipotent characters, denoted
in~\cite{Ward} by
$\xi_k$ for $k\in I$. In~\cite{LusCox}, Lusztig shows that $\omega_{\xi_1}=1$,
$\omega_{\xi_3}=1$, and
$$\{\omega_{\xi_5},\omega_{\xi_{7}}\}=\{\pm i\}\quad
\quad\textrm{and}\quad
\{\omega_{\xi_9},\omega_{\xi_{10}}\}=\left\{\frac{\pm
i-\sqrt{3}}{2}\right\},$$
where $i\in\overline{\Q}_{\ell}$ is a root of $-1$. This work is
completed in~\cite{GM} by Geck-Malle who proved that 
$\{\omega_{\xi_6},\omega_{\xi_{8}}\}=\{\pm i\}$.

The aim of this note is to compute the roots $\omega_{\xi_k}$ for $k\in I$.
Moreover, we will also 
show that the Fourier matrices of the Ree groups of type $G_2$ satisfy
a conjecture of
Digne-Michel~\cite{DMShin} that we recall in~\S\ref{part4}.
These are new results, which
complete 
works of Lusztig~\cite{LusCox} and of Geck-Malle~\cite{GM} for the
Ree groups of type $G_2$.

The paper is organized as follows. In Section~\ref{part1} we fix some
notation and we give preliminary results. In Section~\ref{part2} we
give results on the Shintani descents from $\Galg^{F^2}\semi F$ to
$\Galg^F$ that we need in order to apply the same method as in~\cite{Br3}.
In Section~\ref{part3} we compute the roots of
unity associated as above to the unipotent characters of the Ree
groups. Finally, in Section~\ref{part4} we show that the Fourier
matrices
for the Ree groups defined by Malle and Geck satisfy the Digne-Michel
conjecture.
\section{Notation and preliminary results}\label{part1}

\subsection{Notation}
Let $\Galg$ be a simple algebraic group of type $G_2$ over an
algebraic closure 
$\overline{\F}_3$ of the finite field $\F_3$ with $3$ elements.
We denote by $\Sigma$ the root system of type $G_2$, and by $\Pi=\{a,b\}$
a fundamental system of roots. We choose $a$ for the short root, and $b$
for the long one. We denote by $\Sigma^+=\{a,b,a+b,2a+b,3a+b,3a+2b\}$  
the set of positive roots with
respect to $\Pi$. For $r\in\Sigma$ and $t\in\overline{\F}_3$, there is the corresponding
Chevalley element $x_r(t)\in\Galg$. We recall that $\Galg=\cyc{x_r(t)\ |\
r\in\Sigma, t\in\overline{\F}_3^{\times}}$. 
We set
$$\Ualg=\cyc{x_r(t)\ |\ r\in\Sigma^+,\,t\in\overline{\F}_3}\quad
\textrm{and}\quad\Talg=\cyc{h_r(t)\ |\
r\in\Sigma^+,\,t\in\overline{\F}_3},$$
where $h_r(t)=x_{-r}(t^{-1}-1)x_r(1)x_{-r}(t-1)x_r(-t^{-1})$. 
The subgroup $\Talg$ is a maximal torus of $\Galg$, contained
in the Borel subgroup $\Balg=\Talg\Ualg$ of $\Galg$.
The Weyl group of $\Galg$ is $W=\Nn_{\Galg}(\Talg)/\Talg$.

For every positive integer $m$,
the Frobenius map $F_m$ on $\Galg$ is defined on the Chevalley
generators by setting
$F_m(x_r(t))=x_r(t^{3^m})$. As in~\cite[\S12.4]{Carter1}, we define an
automorphism $\alpha$ of $\Galg$ by setting 
$\alpha(x_{r}(t))=x_{\rho(r)}(t)$ if $r$ is a long root and 
$\alpha(x_{r}(t))=x_{\rho(r)}(t^3)$ if $r$ is short, where $\rho$ is the
unique angle-preserving, and length-changing bijection of $\Sigma$ which
preserves $\Pi$.

Throughout this paper, we fix a positive integer $n$. 
We set $\theta=3^n$ and $q=3\theta^2$. We write $F=\alpha\circ
F_n$. We then have $F^2=F_{2n+1}$.
The fixed-point subgroups $\Galg^F$ and $\Galg^{F^2}$ are the Ree
group $^2G_2(q)$ and the finite Chevalley group $G_2(q)$
respectively.
The subgroups $\Talg$, $\Ualg$ and $\Balg$ are $F$-stable. We notice
that the automorphism of $W$ induced by $F$ has order $2$.

%
%
Moreover, the
Chevalley relations are, for $u,v\in\overline{\F}_3$:
$$\begin{array}{lll}
x_a(t)x_b(u)&=&x_b(u)x_a(t)x_{a+b}(-tu)x_{3a+b}(t^3u)x_{2a+b}(-t^2u)x_{3a+2b}(t^3u^2)\\
x_{a}(t)x_{a+b}(u)&=&x_{a+b}(u)x_a(t)x_{2a+b}(tu)\\
x_{b}(t)x_{3a+b}(u)&=&x_{3a+b}(u)x_b(t)x_{3a+2b}(tu)\\
x_{a+b}(t)x_{3a+b}(u)&=&x_{3a+b}(u)x_{a+b}(t)\\
x_{a+b}(t)x_{2a+b}(u)&=&x_{2a+b}(u)x_{a+b}(t)\\
x_{a+b}(t)x_{3a+2b}(u)&=&x_{3a+2b}(u)x_{a+b}(t)\\
x_{2a+b}(t)x_{3a+b}(u)&=&x_{3a+b}(u)x_{2a+b}(t)\\
x_{2a+b}(t)x_{3a+2b}(u)&=&x_{3a+2b}(u)x_{2a+b}(t)
\end{array}$$

We fix a root
$\alpha_0\in\overline{\F}_3^\times$
of $X^q-X+1$, and we set $\xi=\alpha_0^3-\alpha_0$. 
We have
$$\xi^q=\alpha_0^{3q}-\alpha_0^q=\alpha_0^3-1-\alpha_0+1=\xi.$$
Therefore $\xi\in\F_q$. Moreover, 
$X^3-X-\xi\in\F_q[X]$ is
irreducible over $\F_q$. Otherwise there is a $t\in\F_q$ with
$t^3-t-\xi=0$, implying $(t-\alpha_0)^3=(t-\alpha_0)$. However,
$t\neq\alpha_0$ (because $\alpha_0\not\in\F_q$). Thus
$(t-\alpha_0)^2=1$. It follows that $\alpha_0=t\pm 1\in\F_q$. This is a
contradiction.  

The character table of $\Galg^{F^2}$ was computed by
Enomoto in~\cite{Enomoto}. The description of this table depends on an
element $\xi_0\in\F_q$ such that $X^3-X-\xi_0$ is an irreducible
polynomial
over $\F_q$.
In the following we
choose $\xi_0=\xi$. 

We recall that the
unipotent regular
class $u_{\textrm{reg}}$ of $\Galg$ splits in $3$ classes $A_{51}$,
$A_{52}$ and $A_{53}$ in $\Galg^{F^2}$, with
representatives $x_a(1)x_b(1)$, $x_a(1)x_b(1)x_{3a+b}(\xi)$ and
$x_a(1)x_b(1)x_{3a+b}(-\xi)$ respectively. 
Moreover, $u_{\textrm{reg}}$ also splits in
$3$ classes in $\Galg^F$ whose representatives are not conjugate in
$\Galg^{F^2}$. We denote by $Y_1$, $Y_2$ and $Y_3$ representatives 
with $Y_1\in A_{51}$, $Y_2\in A_{52}$ and $Y_3\in A_{53}$.

The group $\Galg^F$ has $q+8$ conjugacy classes. More
precisely, we recall 
that the system of representatives of 
classes of $\Galg^F$ given in~\cite[4.1]{Br5} is described 
as follows. 

\begin{itemize}
\item The trivial element $1$.
\item The element $J=h_{a+b}(-1)$ which has a centralizer of order
$q(q-1)(q+1)$.
\item The element $X=x_{2a+b}(1)x_{3a+2b}(1)$ which has 
centralizer order $q^3$.
\item The elements $T=x_{a+b}(1)x_{3a+b}(1)$ and $T^{-1}$
whose centralizers have
order $2q^2$.
\item The elements $Y_1$, $Y_2$ and $Y_3$ described as above whose
centralizers have order $3q$.
\item The elements $TJ$ and $T^{-1}J$ whose centralizers have order
$2q$.
\item A family $R$ of $(q-3)/2$ semisimple regular elements with 
centralizer order $q-1$.
\item A family $S$ of $(q-3)/6$ semisimple regular elements with 
centralizer order $q+1$.
\item A family $V$ of $(q-3\theta)/6$ semisimple regular elements with 
centralizer order $q-3\theta+1$.
\item A family $M$ of $(q+3\theta)/2$ semisimple regular elements with 
centralizer order $q+3\theta+1$.
\end{itemize}

In~\cite[4.5]{Br5} we give 
the class fusion between
$\Galg^F$ and $\Galg^{F^2}$.
The character table of $\Galg^F$ with respect 
to this parametrization is given in~\cite[7.2]{Br5}. We notice that
there are some misprints in~\cite[7.2]{Br5} for the values of $\xi_9$
and $\xi_{10}$ on $Y_2$ and $Y_3$. Indeed we have
$$\xi_9(Y_2)=\xi_{10}(Y_3)=\theta(1+i\sqrt{3})/2\quad\textrm{and}\quad
\xi_9(Y_3)=\xi_{10}(Y_2)=\theta(1-i\sqrt{3})/2.$$

A system of representatives of the conjugacy classes of
$\Galg^{F^2}\semi F$ is computed in~\cite[4.2]{Br5}. It is shown 
that the following elements are representatives of the
outer classes of $\Galg^{F^2}\semi F$ (i.e. the classes of
$\Galg^{F^2}\semi F$ lying in the coset $\Galg^{F^2}.F$): 

\begin{itemize}
\item The element $F$, which has centralizer order 
$2q^3(q^2-1)q^2-q+1)$.
\item The element $h_0.F$ with $h_0=h_{a}(-1)$,
 which has  centralizer order 
$2q(q-1)(q+1)$.
\item The element $X.F$, with  
centralizer order $2q^3$.
\item The elements $T.F$ and $T^{-1}.F$, whose centralizers have
order  $6q^2$.
\item The elements $Y_1.F$, $Y_2.F$ and $Y_3.F$, whose
centralizers have order $6q$.
\item The elements $\eta h_0.F$ and $\eta^{-1} h_0.F$ with
$\eta=x_{a+b}(1)x_{3a+b}(-1)$,
whose centralizers have order
$2q$.
\item A family $R'$ of $(q-3)/2$ elements with 
centralizer order $q-1$.
\item A family $S'$ of $(q-3)/6$ elements with
centralizer order $q+1$.
\item A family $V'$ of $(q-3\theta)/6$ elements with 
centralizer order $q-3\theta+1$.
\item A family $M'$ of $(q+3\theta)/2$ elements with 
centralizer order $q+3\theta+1$.
\end{itemize}

Finally, we recall that the character table of $\Galg^{F^2}\semi F$ 
is computed in~\cite[1.1]{Br5}. 


\subsection{Some results}

We will need the following results.

\begin{lemma}\label{conjU}
We set $E=\{t^3-t\ |\ t\in\F_q\}$. Then every element
$x\in\F_q$ can be written uniquely as $x=\eta_x\xi+y_x$ with $y_x\in E$ and
$\eta_x\in\F_3$. 
Moreover let $u\in\Ualg^{F^2}$ be such that
$$u=x_a(1)x_b(1)x_{a+b}(t_{a+b})x_{3a+b}(t_{3a+b})x_{2a+b}(t_{2a+b})
x_{3a+2b}(t_{3a+2b}),$$ for $t_{a+b}$, $t_{3a+b}$, $t_{2a+b}$,
$t_{3a+2b}\in\F_q$. For $u$ as above, we set $p(u)=t_{a+b}+t_{3a+b}$. 
Then $u\in A_{51}$ if $\eta_{p(u)}=0$, $u\in A_{52}$ if $\eta_{p(u)}=1$, and
$u\in A_{53}$ if $\eta_{p(u)}=-1$.
\end{lemma}
\begin{proof}This is a consequence of the Chevalley relations.
\end{proof}

We remark that the map $\F_q\rightarrow\F_3,\,x\mapsto\eta_x$ is an
additive morphism.
We now describe the elements $Y_1$, $Y_2$ and $Y_3$ more precisely.

\begin{lemma}\label{conjRep}
For $u=\pm\xi$, we set
$$\alpha(1)=x_a(1)x_b(1)x_{a+b}(1)x_{2a+b}(1)\ \textrm{and }
\beta(u)=x_{a+b}(u^\theta)x_{3a+b}(u).$$
As previously,
we denote by $\eta_1\in\F_3$ the unique element such that
$1=\eta_1\xi+t^3-t$ for some $t\in\F_q$.
\begin{itemize}
\item If $n\equiv 
1\mod 3$, then $\eta_1=0$. We choose
$Y_1=\alpha(1)$, 
$Y_2=\alpha(1)\beta(-\xi)$ and
$Y_3=\alpha(1)\beta(\xi)$.
\item If $n\equiv 
0\mod 3$, then $\eta_1=-1$. We choose
$Y_1=\alpha(1)\beta(-\xi)$, 
$Y_2=\alpha(1)\beta(\xi)$ and
$Y_3=\alpha(1)$.
\item If $n\equiv 
-1\mod 3$, then $\eta_1=1$. We choose
$Y_1=\alpha(1)\beta(\xi)$, 
$Y_2=\alpha(1)$ and
$Y_3=\alpha(1)\beta(-\xi)$.
\end{itemize}
\end{lemma}

\begin{proof}We have $$p(\alpha(1))=1,\quad
p(\alpha(1)\beta(\xi))=1+\xi+\xi^{\theta}\quad\textrm{and}\quad
p(\alpha(1)\beta(-\xi))=1-\xi-\xi^{\theta}.$$ 
We discuss $\eta_1$. There is an element $t\in\F_q$ such that
$$1=\eta_1\xi+t^3-t.$$
For $0\leq j\leq 2n$, we take the $3^j$-power of the last relation,
and sum all new obtained relations. Thus we obtain
$$2n+1=\eta_1(\xi+\xi^3+\cdots+\xi^{3^{2n}}).$$
However, $\xi=\alpha_0^3-\alpha_0$ with $\alpha_0^q=\alpha_0-1$. Hence
$\xi+\cdots+\xi^{3^{2n}}=\alpha_0^q-\alpha_0=-1$. It follows that
$$2n+1+\eta_1\equiv 0\mod 3.$$ 
Moreover, we remark that
$\xi^\theta=(\xi+\cdots+\xi^{\theta/3})^{3}-(\xi+\cdots+\xi^{\theta/3})
+\xi$. We deduce that $\eta_{\xi^{\theta}}=1$.
\begin{itemize}
\item If $n\equiv 1\mod 3$, then $\eta_1=0$. Therefore
$\eta_{1+\xi+\xi^{\theta}}=-1$ and $\eta_{1-\xi-\xi^{\theta}}=1$.
\item If $n\equiv 0\mod 3$, then $\eta_1=-1$. Therefore
$\eta_{1+\xi+\xi^{\theta}}=1$ and $\eta_{1-\xi-\xi^{\theta}}=0$.
\item If $n\equiv -1\mod 3$, then $\eta_1=1$. Therefore
$\eta_{1+\xi+\xi^{\theta}}=0$ and $\eta_{1-\xi-\xi^{\theta}}=-1$.
\end{itemize}
The result follows.
\end{proof}

\begin{lemma}\label{systeme}Let $\alpha$ and $\beta$ 
be two elements of $\F_q$. We consider
the following system of equations (S)
$$\left\{\begin{array}{lcl}
y^\theta-x&=&1\\
x^{3\theta}-y&=&1\\
t^{\theta}-z+1+x^{3\theta+1}&=&\alpha\\
z^{3\theta}-t-1-x^{3\theta+3}&=&\beta
\end{array}\right.$$
If $x_0$ is a root in $\overline{\F}_3$ of $X^q-X+1$ and $t_0$
is a root in $\overline{\F}_3$ of
$X^q-X-x_0^{3\theta}-\beta-\alpha^{3\theta}$, then the tuple
$(x_0,x_0^{3\theta}-1,t_0,x_0^{3\theta+1}+t_0^\theta+1-\alpha)$ is a
solution of (S) in $\overline{\F}_3$.
\end{lemma}

\begin{proof}If $y_0=x_0^{3\theta}-1$, then
$y_0^\theta=x_0^q-1=x_0+1$. If $z_0=x_0^{3\theta+1}+t_0^\theta+1-\alpha$,
then $$\begin{array}{lcl}
z_0^{3\theta}&=&x_0^{3q}x_0^{3\theta}+t_0^q+1-\alpha^{3\theta}\\
&=&(x_0^3-1)x_0^{3\theta}+t_0+x_0^{3\theta}+b+\alpha^{3\theta}+1-\alpha^{3\theta}\\
&=&x_0^{3\theta+3}+t_0+\beta+1
\end{array}$$
The result follows.
\end{proof}
\begin{remark}\label{remarque1} 
For the solution of the system $(S)$ given in
Lemma~\ref{systeme}, we remark that 
we can choose $x_0$ independently of $\alpha$ and
$\beta$.
\end{remark}
\section{Shintani descents}\label{part2}

\subsection{Definition}

A reference for this section is for example~\cite{DMShin}. 
We keep the notation as
above. We denote by
$L_F:\Galg\rightarrow\Galg,\,x\mapsto x^{-1}F(x)$ the Lang map
associated to $F$.
Since $\Galg$ is a simple algebraic group, it follows that $\Galg$ is
connected. Hence $L_F$ and $L_{F^2}$ are surjective. Moreover, 
$L_{F}(x)\in\Galg^{F^2}$ for some $x\in\Galg$ if and only if
$L_{F^2}(x^{-1})\in\Galg^F$, and the correspondence
$$L_{F^2}(x^{-1})\in\Galg^{F}\longleftrightarrow
L_F(x).F\in\Galg^{F^2}\semi F\quad\textrm{for }x\in\Galg,$$
induces a bijection from the outer classes of the group
$\Galg^{F^2}\semi F$ to
the 
classes of $\Galg^F$. We denote this correspondence by $N_{F/F^2}$.
Furthermore we have~\cite[\S I.7]{DMShin}
\begin{equation}\label{eqcent}
|\Cen_{\Galg^{F^2}\semi F}\left( L_{F}(x).F
\right)|=2|\Cen_{\Galg^F}(L_{F^2}(x^{-1})|\quad
\textrm{for }
x\in\Galg.
\end{equation}

For every class function $\psi$ on $\Galg^{F^2}\semi F$, we define the
Shintani $\Sh_{F^2/F}(\psi)$ of $\psi$ by
$\Sh_{F^2/F}(\psi)=\psi\circ N_{F/F^2}$.

\subsection{Shintani correspondence from $\Galg^{F}$ to
$\Galg^{F^2}\semi F$ in type
$G_2$}

\begin{lemma}\label{shinT}
We keep the notation as above. We set $T=\beta(1)$.
We have
$$N_{F/F^2}([T]_{\Galg^F})=[T.F]_{\Galg^{F^2}\semi F}\quad
\textrm{and}\
N_{F/F^2}([T^{-1}]_{\Galg^F})=[T^{-1}.F]_{\Galg^{F^2}\semi
F}.$$
Here $[x]_G$ denotes the conjugacy class of $x$ in $G$.
\end{lemma}

\begin{proof} We set
$x=x_{a+b}(\alpha_0)x_{3a+b}(\alpha_0^{3\theta}-1)$. Then we have
$$\begin{array}{lcl}
L_F(x)&=&x_{a+b}(
(\alpha_0^{3\theta}-1)^{\theta}-\alpha_0)x_{3a+b}(\alpha_0^{3\theta}
-(\alpha_0^{3\theta}-1))\\
&=&x_{a+b}(
\alpha_0-2-\alpha_0)x_{3a+b}(1)\\
&=&\beta(1)
\end{array}
$$ and 
$L_{F^2}(x^{-1})=x_{a+b}(\alpha_0-\alpha_0^q)
x_{3a+b}(\alpha_0-\alpha_0^q)=\beta(1)$.
\end{proof}

\begin{lemma}\label{shinY}
We keep the notation as above. Then we have
$$\begin{array}{lcl}
N_{F/F^2}([Y_1]_{\Galg^F})&=&[Y_3.F]_{\Galg^{F^2}\semi F}\\
N_{F/F^2}([Y_2]_{\Galg^F})&=&[Y_1.F]_{\Galg^{F^2}\semi F}\\
N_{F/F^2}([Y_3]_{\Galg^F})&=&[Y_2.F]_{\Galg^{F^2}\semi F}.
\end{array}$$
%
\end{lemma}

\begin{proof}Let
$u=x_a(1)x_b(1)x_{a+b}(\alpha)x_{3a+b}(\beta)y\in\Ualg^{F^2}$ with $y\in
\operatorname{Z}(\Ualg^{F^2})$.
Since $\Ualg$ is connected, there is $x\in\Ualg$ such that $L_F(x)=u$. 
Then there are $z\in\operatorname{Z}(\Ualg^{F^2})$ and $t_a$, $t_b$,
$t_{a+b}$ and $t_{3a+b}$ in $\overline{\F}_3$ such that
$$x=x_a(t_a)x_b(t_b)x_{a+b}(t_{a+b})x_{3a+b}(t_{3a+b})z.$$
Using the Chevalley relations, we have
$$\begin{array}{lcl}L_F(x)&=&x_a(t_b^\theta-t_a)x_b(t_a^{3\theta}-t_b)
x_{a+b}(t_at_b-t_{a+b}+ t_b^\theta
t_a^{3\theta}+t_{3a+b}^{\theta}-t_b^{\theta+1})\\
&& x_{3a+b}(t_{a+b}^{3\theta}-t_b^{3\theta}t_a^{3\theta}
-t_a^3t_b-t_{3a+b}+t_b^{3\theta+1})z',
\end{array}
$$
for some $z'\in\operatorname{Z}(\Ualg^{F^2})$.
However, using the uniqueness of the Chevalley decomposition we have the
system of equations $(S')$ 
$$\left\{
\begin{array}{lcl}
t_b^\theta-t_a&=&1\\
t_a^{3\theta}-t_b&=&1\\
t_at_b-t_{a+b}+ t_b^\theta
t_a^{3\theta}+t_{3a+b}^{\theta}-t_b^{\theta+1}&=&\alpha\\
t_{a+b}^{3\theta}-t_b^{3\theta}t_a^{3\theta}
-t_a^3t_b-t_{3a+b}+t_b^{3\theta+1}&=&\beta
\end{array}
\right.$$
We deduce from the two first equations that $t_a^q-t_a=-1$ and
$t_b^q-t_b=-1$.

Moreover, there is $z''\in\operatorname{Z}(\Ualg^{F^2})$ such that
$$
\begin{array}{lcl}L_{F^{2}}(x^{-1})&=&
x_a(t_a-t_a^q)x_b(t_b-t_b^q)
x_{a+b}(t_a^q(t_b^q-t_b) + t_{a+b} - t_{a+b}^{q})\\
&& x_{3a+b}(t_{3a+b}-t_{3a+b}^q-t_a^{3q}(t_b-t_b^q))z''\\
&=&x_a(1)x_b(1)x_{a+b}(t_{a+b}-t_{a+b}^q
-t_a^q)x_{3a+b}(t_{3a+b}-t_{3a+b}^q+t_a^{3q})z''.
\end{array}
$$
Hence using Lemma~\ref{conjRep} in order to find the $\Galg^F$-class
of $L_{F^{2}}(x^{-1})\in\Galg^F$, it is sufficient to find the 
$\Galg^{F^2}$-class
of $L_{F^{2}}(x^{-1})$. However using Lemma~\ref{conjU}, we have to
evaluate
$\eta_{p(L_{F^2}(x^{-1}))}$ where $p$ is defined as in Lemma~\ref{conjU}.
We have 
$$\begin{array}{lcl}
p(L_{F^2}(x^{-1}))&=&t_{3a+b}-t_{3a+b}^q+t_a^{3q}+t_{a+b}-t_{a+b}^q
-t_a^q\\
&=&t_{3a+b}-t_{3a+b}^q+t_{a+b}-t_{a+b}^q+t_a^3-t_a.
\end{array}$$
Using the equations of the system $(S')$, we deduce that 
$$\left\{
\begin{array}{lcl}
t_a^{3\theta+1}+t_{3a+b}^{\theta}-t_{a+b}+1&=&\alpha\\
-t_a^{3\theta+3}+t_{a+b}^{3\theta}-t_{3a+b}-1&=&\beta
\end{array}\right.
$$
It follows that
$$\left\{
\begin{array}{lcl}
t_a^{(3\theta+1)3\theta}-t_a^{3\theta+3}+t_{3a+b}^q-t_{3a+b}&=&\alpha^{3\theta}+\beta\\
t_a^{3\theta+1}-t_a^{(3\theta+3)\theta}+t_{a+b}^q-t_{a+b}&=&\alpha+\beta^{\theta}
\end{array}\right.
$$
Adding these two equations, we obtain
$$t_{a+b}^q-t_{a+b}+t_{3a+b}^q-t_{3a+b}=\alpha+\alpha^{3\theta}+
\beta+\beta^{\theta}.$$
Hence we have
$$p(L_{F^2}(x^{-1}))=-(\alpha+\alpha^{3\theta}+
\beta+\beta^{\theta})+t_a^3-t_a.$$
Moreover, we remark that the system $(S')$ is equivalent to the system
$(S)$. We use for $(t_a,t_b,t_{a+b},t_{3a+b})$ the solution described 
in Lemma~\ref{systeme}, choosing
$t_a=\alpha_0$, which is possible as we have seen in
Remark~\ref{remarque1}. Thus
$$p(L_{F^2}(x^{-1}))=-(\alpha+\alpha^{3\theta}+
\beta+\beta^{\theta})+\xi.$$
We suppose now that $n\equiv 1\mod 3$. Using Lemma~\ref{conjRep}, we
have 
$$
\begin{array}{lcl}
Y_1&=&x_a(1)x_b(1)x_{a+b}(1)x_{2a+b}(1)\\
Y_2&=&x_a(1)x_b(1)x_{a+b}(1-\xi^\theta)x_{3a+b}(-\xi)x_{2a+b}(1)\\
Y_3&=&x_a(1)x_b(1)x_{a+b}(1+\xi^\theta)x_{3a+b}(\xi)x_{2a+b}(1)
\end{array}$$
For $(\alpha,\beta)=(1,0)$,
$(\alpha,\beta)=(1-\xi^{\theta},-\xi)$ and 
$(\alpha,\beta)=(1+\xi^{\theta},\xi)$ respectively, we deduce that
$\eta_{p(L_{F^2}(x^{-1}))}$ is equal to
$1$, $-1$, and~$0$ respectively, because $\eta_1=0$. Thus 
$$\begin{array}{lcl}
N_{F/F^2}([Y_1]_{\Galg^F})&=&[Y_3.F]_{\Galg^{F^2}\semi F}\\
N_{F/F^2}([Y_2]_{\Galg^F})&=&[Y_1.F]_{\Galg^{F^2}\semi F}\\
N_{F/F^2}([Y_3]_{\Galg^F})&=&[Y_2.F]_{\Galg^{F^2}\semi F}
\end{array}$$
We proceed similarly when $n\equiv 0\mod 3$ and $n\equiv -1\mod 3$.
The result follows.\end{proof}

\begin{proposition}\label{corrshin}
We keep the same notation as above. Then we have
$$\begin{array}{lcl}
N_{F/F^2}([1]_{\Galg^F})&=&[F]_{\Galg^{F^2}\semi F}\\
N_{F/F^2}([X]_{\Galg^F})&=&[X.F]_{\Galg^{F^2}\semi F}\\
N_{F/F^2}([J]_{\Galg^F})&=&[h_0.F]_{\Galg^{F^2}\semi F}\\
N_{F/F^2}([T]_{\Galg^F})&=&[T.F]_{\Galg^{F^2}\semi F}\\
N_{F/F^2}([T^{-1}]_{\Galg^F})&=&[T^{-1}.F]_{\Galg^{F^2}\semi F}\\
N_{F/F^2}([Y_1]_{\Galg^F})&=&[Y_3.F]_{\Galg^{F^2}\semi F}\\
N_{F/F^2}([Y_2]_{\Galg^F})&=&[Y_1.F]_{\Galg^{F^2}\semi F}\\
N_{F/F^2}([Y_3]_{\Galg^F})&=&[Y_2.F]_{\Galg^{F^2}\semi F}\\
\{N_{F/F^2}([JT]_{\Galg^F}),N_{F/F^2}([JT^{-1}]_{\Galg^F})\}&=&
\{[\eta 
h_0.F]_{\Galg^{F^2}\semi F},[\eta^{-1}h_0.F]_{\Galg^{F^2}\semi F}\}\\
\{N_{F/F^2}([t]_{\Galg^F})\ |\ t\in R\}&=&\{[t]_{\Galg^{F^2}\semi F}\ |\
t\in R'\}\\
\{N_{F/F^2}([t]_{\Galg^F})\ |\ t\in S\}&=&\{[t]_{\Galg^{F^2}\semi F}\ |\
t\in S'\}\\
\{N_{F/F^2}([t]_{\Galg^F})\ |\ t\in V\}&=&\{[t]_{\Galg^{F^2}\semi F}\ |\
t\in V'\}\\
\{N_{F/F^2}([t]_{\Galg^F})\ |\ t\in M\}&=&\{[t]_{\Galg^{F^2}\semi F}\ |\
t\in M'\}\\
\end{array}$$

\end{proposition}
\begin{proof}To prove this result, we essentially use
Relation~(\ref{eqcent}) comparing the orders of centralizers of the
representatives  of classes of $\Galg^F$, and of the representatives 
of the outer classes of $\Galg^{F^2}\semi F$
recalled
in \S\ref{part1}. Moreover, for the classes $[T]_{\Galg^F}$, 
$[T^{-1}]_{\Galg^F}$, and  $[Y_1]_{\Galg^F}$,  $[Y_2]_{\Galg^F}$,
$[Y_3]_{\Galg^F}$, we use Lemma~\ref{shinT} and Lemma~\ref{shinY}
respectively.
\end{proof}

\subsection{Shintani descents of the unipotent characters in type
$G_2$}

The $F$-stable unipotent characters of $\Galg^{F^2}$ are denoted by
$1_{\Galg^{F^2}}$, $\theta_1$, $\theta_2$, $\theta_5$, $\theta_{10}$,
$\theta_{11}$, $\theta_{12}[1]$ and $\theta_{12}[-1]$
in~\cite{Enomoto}. Their degrees are $1$,
$\frac{1}{6}q(q+1)^2(q^2+q+1)$, $\frac{1}{2}q(q+1)(q^3+1)$, $q^6$,
$\frac{1}{6}q(q-1)^2(q^2-q+1)$, $\frac{1}{2}q(q-1)(q^3-1)$,
$\frac{1}{3}q(q^2-1)^2$ and $\frac{1}{3}q(q^2-1)^2$ respectively.
These characters extends to $\Galg^{F^2}\semi F$. Let $\chi$ be such a
character. If $\chi$ has an extension such that its value on $F$ is
not zero, then we denote by $\widetilde{\chi}$ the extension of $\chi$
such that $\widetilde{\chi}(F)>0$. 
In \cite[5.6]{Br5} we have seen that this is always the case except
for $\theta_2$ and $\theta_{10}$. Then we choose for
$\widetilde{\theta}_2$ and $\widetilde{\theta}_{10}$ the extensions
such that
$$\widetilde{\theta}_2(\eta
h_0.F)=\sqrt{q}\quad\textrm{and}\quad\widetilde{\theta}_{10}(T.F)=
\theta^2\sqrt{3}i.$$ 
Moreover, there is a misprint in~\cite[5.6]{Br5}. Indeed,
we have $\widetilde{\theta}_{10}(Y_2)=-\theta\sqrt{3}i$ and
$\widetilde{\theta}_{10}(Y_3)=\theta\sqrt{3}i$.

\begin{proposition}\label{desShin}
We keep the same notation as above. Then we have
$$\begin{array}{lcl}
\Sh_{F^2/F}(1_{\Galg^{F^2}\semi F})&=&1_{\Galg^F}\\
\Sh_{F^2/F}(\widetilde{\theta}_1)&=&\dfrac{\sqrt{3}}{6}\left(i\xi_5
+i\xi_6-i\xi_7-i\xi_8+(\sqrt{3}-i)\xi_9+(\sqrt{3}+i)\xi_{10}\right)\\
\Sh_{F^2/F}(\widetilde{\theta}_2)&=&\pm\dfrac{1}{2}(-i\xi_5+i\xi_6+i\xi_7-i\xi_8)\\
\Sh_{F^2/F}(\widetilde{\theta}_{5})&=&\xi_3\\
\Sh_{F^2/F}(\widetilde{\theta}_{10})&=&\dfrac{\sqrt{3}}{6}\left(i\xi_5
+i\xi_6+i\xi_7+i\xi_8+(\sqrt{3}-i)\xi_9-(i+\sqrt{3})\xi_{10}\right)\\
\Sh_{F^2/F}(\widetilde{\theta}_{11})&=&\dfrac{1}{2}(\xi_5-\xi_6+\xi_7
-\xi_8)\\
\Sh_{F^2/F}(\widetilde{\theta}_{12}[1])&=&\dfrac{\sqrt{3}}{6}\left(
 (\sqrt{3}-i)\xi_5+ (\sqrt{3}-i)\xi_6+(\sqrt{3}+i)\xi_9\right)\\
\Sh_{F^2/F}(\widetilde{\theta}_{12}[-1])&=&\dfrac{\sqrt{3}}{6}\left(
 (\sqrt{3}+i)\xi_7+ (\sqrt{3}+i)\xi_8+(\sqrt{3}-i)\xi_{10}\right)\\
\end{array}$$
\end{proposition}
\begin{proof}
We use Proposition~\ref{corrshin} and \cite[5.6]{Br5} to obtain the values 
of the Shintani descents as class functions on $\Galg^F$. 
Using
the character table of $\Galg^F$ in~\cite{Ward}, we decompose them in
the basis of irreducible characters of $\Galg^F$. The result follows.
\end{proof}

\section{Eigenvalues of the Frobenius for the Ree groups of type $G_2$}
\label{part3}

As an application of Proposition~\ref{desShin} we compute in this section
the roots of unity associated to the unipotent characters of the Ree groups
of type $G_2$, as explained
in Section~\ref{intro}.

\begin{theorem}We keep the same notation as above. The roots of unity
associated to the unipotent characters of $^2G_2(q)$ are
$$\renewcommand{\arraystretch}{1.4}\begin{array}{l|l|l|l|l|l|l|c|c}
\chi&\xi_1&\xi_3&\xi_5&\xi_6&\xi_7&\xi_8&\xi_9&\xi_{10}\\
\hline
\omega_{\chi}&1&1&-i&-i&i&i&\frac{-\sqrt{3}+i}{2}&\frac{-\sqrt{3}-i}{2}
\end{array}$$

\end{theorem}
\begin{proof}We first recall a result of Digne-Michel. 
Let $\rho$ be an irreducible character of $W$. Using the Harish-Chandra theory, we can associate to
$\rho$ an irreducible character $\chi_{\rho}$ of the principal series
of $\Galg^{F^2}$, that is the set of irreducible constituents of
$$\Phi=\Ind_{\Balg^{F^2}}^{\Galg^{F^2}}(1_{\Balg^{F^2}}).$$
The group $\cyc F$ acts on $\Irr(W)$. If $\rho$ is
$F$-stable, then $\rho$ extends to $W\semi F$, and has exactly $2$
extensions denoted by $\widetilde{\rho}$ and
$\varepsilon\widetilde{\rho}$, where $\varepsilon$ is the non-trivial
character of $W\semi F$ which has $W$ is its kernel. 
However Malle shows in~\cite[1.5]{MalleHarish} that the irreducible characters of $W\semi F$ are 
in 1-1 correspondence with the constituents of
$\widetilde{\Phi}=\Ind_{\Galg^{F^2}}^{\Galg^{F^2}\semi F}(\Phi)$. In
particular, if $\rho$ is $F$-stable, then so is $\chi_{\rho}$. Hence, the
characters $\chi_{\widetilde{\rho}}$ and
$\chi_{\widetilde{\varepsilon\rho}}$ of  
$\widetilde{\Phi}$ corresponding to $\widetilde{\rho}$ and
$\varepsilon\widetilde{\rho}$ respectively, are the two extensions of
$\chi_{\rho}$ to $\Galg^{F^2}\semi F$. Moreover, we recall that the
almost character of $\Galg^F$ corresponding to $\widetilde{\rho}$ is
defined by
$$\mathcal{R}_{\widetilde{\rho}}=\sum_{w\in
W}\widetilde{\rho}(w.F)R_w.$$
The main theorem~\cite[2.3]{DMShin} asserts that
\begin{equation}
\Sh_{F^2/F}(\chi_{\widetilde\rho})=\sum_{V\in
\mathcal{U}(\Galg^F}\cyc{\mathcal{R}_{\widetilde{\rho}},V}_{\Galg^F}
\,\omega_V\, V.
\label{DMtheo}
\end{equation}
Furthermore the almost characters of the Ree groups are computed by
Geck and Malle in~\cite[2.2]{GM}. More precisely, the $F$-stable
characters of $W$ are $1_W$, $\operatorname{sgn}$, and the two
characters of degree $2$ of $W$, denoted by $2_1$ and $2_2$. The
character $2_1$ is chosen such that it takes the value $-2$ on the
Coxeter element of $W$. Then for the
extensions of these characters chosen in~\cite{GM}, we have
$$\begin{array}{lcl}
\mathcal{R}_{1_{W\semi F}}&=&1_{\Galg^F}\\ 
\mathcal{R}_{\widetilde{\operatorname{sgn}}}&=&\xi_3\\ 
\mathcal{R}_{\widetilde{2}_1}&=&
\dfrac{\sqrt{3}}{6}(\xi_5+\xi_6+\xi_7+\xi_8+2\xi_9+2\xi_{10})\\ 
\mathcal{R}_{\widetilde{2}_2}&=&\dfrac{1}{2}(\xi_5-\xi_6+\xi_7-\xi_8) 
\end{array}$$
Moreover, using~\cite[p.112,p.150]{Carter2}, we deduce that
$$\Ind_{\Balg^{F^2}}^{\Galg^{F^2}}(1_{\Balg^F})=
1_{\Galg^{F^2}}+\theta_5+\theta_3+\theta_4+2\theta_1+2\theta_2,$$
and $1_{\Galg^{F^2}}=\chi_{1_W}$, $\theta_5=\chi_{\operatorname{sgn}}$, 
$\theta_1=\chi_{2_1}$ and $\theta_2=\chi_{2_2}$. Hence we have
$$\chi_{\widetilde{2}_1}\in\{\widetilde{\theta}_1,\varepsilon\widetilde{\theta}_1\},$$where
$\varepsilon$ denotes now the non trivial character of
$\Galg^{F^2}\semi F$ containing $\Galg^{F^2}$ in its kernel.
Therefore, using Relation~(\ref{DMtheo}) we deduce that
$$\pm\Sh_{F^2/F}(\widetilde{\theta}_1)=\frac{\sqrt{3}}{6}\left(
\omega_{\xi_5}\xi_5+\omega_{\xi_6}\xi_6+\omega_{\xi_7}\xi_7+
\omega_{\xi_8}\xi_8+
2\omega_{\xi_9}\xi_9+2\omega_{\xi_{10}}\xi_{10}\right).$$
Hence we deduce using Proposition~\ref{desShin} that
$$\omega_{\xi_9}=\pm\sqrt{3}\cyc{\Sh_{F^2/F}(\theta_1),\xi_9}_{\Galg^F}
=\pm\frac{\sqrt{3}-i}{2}.$$
However using the result of Lusztig~\cite{LusCox}, we know that
$\omega_{\xi_9}=(\pm i-\sqrt{3})/2$. We then deduce that
$\chi_{\widetilde{2}_1}=\varepsilon\widetilde{\theta}_1$ and that 
$$\omega_{\xi_9}=\frac{-\sqrt{3}+i}{2}.$$
We immediately obtain the other roots using Proposition~\ref{desShin}.
\end{proof}

\begin{remark}To determine the roots of unity of the unipotent
characters of the Ree groups of type $G_2$, the preceding proof shows
that we only need to know the Shintani descent of
$\widetilde{\theta}_1$. 
\end{remark}

\section{Conjecture of Digne-Michel for the Ree groups of type $G_2$}
\label{part4}

In~\cite{DMShin} Digne and Michel state a conjecture on the
decomposition in irreducible constituents of the Shintani descents. We
recall this conjecture in the special case that $
F\in\operatorname{Aut}(W)$ has order $2$. If
$\chi$ is an $F$-stable irreducible character of $\Galg^{F^2}$, and if
$\widetilde{\chi}$ denotes an extension of $\chi$ to $\Galg^{F^2}\semi
F$, then it is conjectured that:
\begin{itemize}
\item The irreducible constituents of $\Sh_{F^2/F}(\widetilde{\chi})$
are unipotent.
\item Up to a normalization, the coefficients of $\Sh_{F^2/F}$ in the
basis~$\Irr(\Galg^F)$ only
depend on the coefficients of the Fourier matrix, and on the roots of
unity attached to the unipotent characters of $\Galg^F$ as above.
More precisely, 
there is a root of unity $u$ such that
$$\pm
u\Sh_{F^2/F}(\widetilde{\chi})=\sum_{V\in\mathcal{U}(\Galg^F)}f_V\omega_V\,V,$$
where $(f_V,\ V\in\mathcal{U}(\Galg^F))$ is, up to a sign, a row of the Fourier matrix.
\end{itemize}

\begin{theorem}\label{conDM}
The conjecture of Digne-Michel on the decomposition of
the Shintani descents of unipotent characters holds in type $G_2$ for
the Frobenius map $F$ that defines the Ree group $^2G_2(q)$.
\end{theorem}

\begin{proof}We set $u_{12[1]}=\frac{1}{2}(-1+\sqrt{3}i)$ and
$u_{12[-1]}=\frac{1}{2}(1+\sqrt{3}i)$. We remark that
$$\begin{array}{ll}
u_{12[1]}(\sqrt{3}-i)=2i,& u_{12[1]}(\sqrt{3}+i)=-\sqrt{3}+i,\\
u_{12[-1]}(\sqrt{3}+i)=2i,& u_{12[-1]}(\sqrt{3}-i)=\sqrt{3}+i
\end{array}$$
Therefore, using Proposition~\ref{desShin} we have
$$
\begin{array}{lcl}
u_{12[1]}\Sh_{F^2/F}(\widetilde{\theta}_{12[1]})&=&
\displaystyle{\frac{\sqrt{3}}{6}
\left(2i\xi_5+2i\xi_6+(i-\sqrt{3})\xi_9\right)}\\
&=&
\displaystyle{\frac{\sqrt{3}}{6}
\left(-2(-i)\xi_5+-2(-i)\xi_6+2\frac{1}{2}(i-\sqrt{3})\xi_9\right).}
\end{array}$$
Similarly, we obtain
$$u_{12[-1]}\Sh_{F^2/F}(\widetilde{\theta}_{12[-1]})=\frac{\sqrt{3}}{6}
\left(2i\xi_7+2i\xi_8-\frac{1}{2}(-i-\sqrt{3})\xi_{10}\right).$$
Moreover, we have
$$\begin{array}{lcl}
\Sh_{F^2/F}(\widetilde{\theta}_2)&=&\pm\dfrac{\sqrt{3}}{6}( 
\sqrt{3}(-i)\xi_5-\sqrt{3}(-i)\xi_6+
\sqrt{3}i\xi_7-\sqrt{3}i\xi_8),\\
\Sh_{F^2/F}(\widetilde{\theta}_{10})&=&\dfrac{\sqrt{3}}{6}\left(-(-i)\xi_5
+-(-i)\xi_6+i\xi_7+i\xi_8-2\dfrac{1}{2}(-\sqrt{3}+i)\xi_9\right.\\
&&\left.+2\dfrac{1}{2}(-i-\sqrt{3})\xi_{10}\right),\\

i\Sh_{F^2/F}(\widetilde{\theta}_{11})&=&\dfrac{\sqrt{3}}{6}
\left(-\sqrt{3}(-i)\xi_5+\sqrt{3}(-i)\xi_6
+\sqrt{3}i\xi_7-\sqrt{3}i\xi_8\right).
\end{array}$$
We set $u_{1}=u_2=u_{10}=1$ and $u_{11}=i$. If we compare the
coefficients in the preceding computations with the coefficients of the
Fourier matrix of Geck-Malle recalled in~\S\ref{intro}, we obtain up to
a sign a row of the Fourier matrix.
Therefore the conjecture holds.
\end{proof}

\begin{remark}\label{rem}
We now discuss the roots $u_i$ for $i\in \{1,2,10,11,12[1],12[-1]\}$ appearing
in the proof of Theorem~\ref{conDM}. 

For $g\in\Galg^F$, Lang's theorem
says that there is $x\in\Galg$ such that $g=x^{-1}F(x)$. Therefore we
have $xgx^{-1}\in\Galg^F$. Moreover, if $x'\in\Galg$ is such that
$x'^{-1}F(x')=g$, then $x'$ lies in the coset $\Galg^F.x$, hence the
$\Galg^F$-class of $xgx^{-1}$ does not depend on the choice of $x$.
For a class function $f$ on $\Galg^F$, we can then define the
Asai-twisting operator $t^*$ by
$$t^*(f)(g)=f(xgx^{-1})\quad\textrm{for }g\in\Galg^F,\ \textrm{and }x\in\Galg\
\textrm{such that }x^{-1}F(x)=g.$$

For a pair $(g,\psi)$ with $g\in\Galg^F$ and $\psi$ an $F$-stable irreducible 
character of the component group 
$A(g)=\Cen_{\Galg}(g)/\Cen_{\Galg}(g)^{\circ}$, we 
can associate a class
function~$\Psi_{(g,\psi)}$ which depends on the choice of an extension
of $\psi$ to $A(g)\semi F$. Therefore $\Psi_{(g,\psi)}$ is an 
eigenvector of $t^*$, and the
corresponding eigenvalue $\lambda_{(g,\psi)}$ is equal to
$\psi(\overline{g})/\psi(1)$, where
$\overline{g}$ denotes the image of $g$ in $A(g)$. 

Let $(c_i,\ i\in I)$ be a row of a Fourier matrix. Then the
construction of the Fourier matrices given by Geck-Malle show that
there is a pair
$(g,\psi)$ as above, such that
$$\Psi_{(g,\psi)}=\sum_{i\in I}c_i\xi_i.$$

To $\chi$ an $F$-stable unipotent character of $\Galg^{F^2}$, we
associate the pair $(g,\psi)$ attached as above to the row of the
Fourier matrix corresponding to $u_{\chi}\Sh_{F^2/F}(\widetilde{\chi})$.
We have
$$\begin{array}{l|c}
\chi&\lambda_{(g,\psi)}\\
\hline
\theta_{1}&1\\
\theta_{2}&1\\
\theta_{10}&1\\
\theta_{11}&-1\\
\theta_{12[1]}&\frac{1}{2}(-1-\sqrt{3}i)\\
\theta_{12[-1]}&\frac{1}{2}(-1+\sqrt{3}i)
\end{array}$$

We observe that the element $u_{\chi}$ chosen in the proof of
Theorem~\ref{conDM} is a root of the polynomial
$$X^2-\lambda_{(g,\psi)}.$$
Moreover, we remark that we can choose for $u_{\chi}$ an arbitrary
root of this polynomial,
because a row of a Fourier matrix is defined up to
a sign.

Finally, we notice that in the situation of type $B_2$
with $F$ the Frobenius map defining the Suzuki groups, this observation
also holds~\cite[4.2]{Br3}. 
\end{remark}
\begin{remark}
Let $\Galg$ be a simple algebraic group of type $F_4$ over
$\overline{\F}_2$, and let $F$ be the Frobenius map on $\Galg$ that 
defines the
Ree group $^2F_4(2^{2n+1})$.
In this situation, the character table of $\Galg^{F^2}\semi F$ is
actually unknown.
Suppose that
\begin{itemize}
\item The conjecture of Digne-Michel holds.
\item We know how to associate to every unipotent character of $^2F_4(q)$
its
root of unity as above.
\item The observation of Remark~\ref{rem} holds.
\end{itemize}
Then, using~\cite[5.4(c)]{GM} and \cite{MaF4}, we can describe 
the values of the
unipotent characters of $\Galg^{F^2}\semi F$ 
on the coset $\Galg^{F^2}.F$.
\end{remark}

\noindent \textit{Acknowledgments.}\quad I wish to thank Gunter Malle for a
helpful and motivating discussion on this work.

%
%
 \bibliographystyle{plain}
 \bibliography{references}
%
%
%
\end{document}